\documentclass{ifacconf}

\usepackage{graphicx}      
\usepackage{natbib}        
\usepackage{amssymb}
\usepackage{amsmath}
\usepackage{dsfont}
\usepackage{euscript}
\usepackage{color}
\usepackage{enumerate}

\newtheorem{assumption}{Assumption}

\newcommand{\R}{\mathbb{R}}

\newcommand{\gam}{\gamma}

\newcommand{\eps}{\epsilon}

\newcommand{\m}{\mathcal}
\newcommand{\mb}{\mathbf}

\begin{document}

\begin{frontmatter}
\title{On infinite dimensional linear programming approach to stochastic control\thanksref{footnoteinfo}}
\thanks[footnoteinfo]{This research is  partially supported by M. Kamgarpour's European Union ERC Starting Grant, CONENE and by T. Summers' the US National Science Foundation under grant CNS-1566127.}
\author[First]{Maryam Kamgarpour}
\author[Second]{Tyler Summers}
\address[First]{Automatic Control Laboratory, ETH Z\"urich, Switzerland (maryamk@ethz.ch)}
\address[Second]{Mechanical Engineering, University of Texas at Dallas (tyler.summers@utdallas.edu)}


\begin{abstract}
We consider the infinite dimensional linear programming (inf-LP) approach for solving stochastic control problems. The inf-LP corresponding to problems with uncountable state and input spaces  is in general computationally intractable. By focusing on linear systems with quadratic cost (LQG), we establish a connection between this approach and the well-known Riccati LMIs. In particular, we show that the semidefinite programs known for the LQG problem  can be derived from the pair of primal and dual inf-LPs. Furthermore, we establish a connection between multi-objective and chance constraint criteria and the inf-LP formulation. 

\end{abstract}


\end{frontmatter}

\section{Introduction}
Optimal control of discrete time stochastic systems can be addressed via the dynamic programming (DP) \citep{bellman1957dp} principle of optimality. For an infinite horizon average or discounted cost problem, the optimal cost function and control policy can be computed as the fixed point of the so-called dynamic programming operator. In general, computing this fixed point is challenging and thus, several approximate approaches based on the DP principle of optimality have been developed. 

An alternative approach to solving stochastic control problems is linear programming (LP) \citep{puterman2009markov,hernandez1996discrete}. If the control and input spaces are uncountable, the corresponding  LP is infinite dimensional (inf-LP). In the primal form of this LP, the optimization variable is the \textit{occupation measure}, which measures infinite horizon occupancy of state and inputs in each Borel subset of the product state input space. An optimal policy may be derived from the optimal occupation measure, while the optimal value function is the optimizer of the dual of this LP. 

In addition to providing an elegant alternative formulation of the optimality conditions for a stochastic control solution, in the LP approach constraints have a natural interpretation. By properly constraining the occupation measure, one can ensure probabilistic constraints on the state trajectory or can ensure bounds on multiple objectives. Such formulations of constrained stochastic control were considered in  \citep{borkar1994ergodic, feinberg1996constrained, altman1999constrained,hernandez2000constrained, hernandez2003constrained}.

The inf-LP formulation is in general computationally intractable. For problems with polynomial data, this inf-LP can be approximated via a sequence of semidefinite programs (SDPs) \citep{savorgnan2009discrete,summers2013approximate}. These recent works are among the few that  explore the inf-LP approach for computation of optimal value function and policies in a stochastic control problem. 

The abstract inf-LP work has not attempted to establish clear connections with the well known, computationally tractable Linear Matrix Inequality (LMI)  formulations of optimal control. In particular, for a stochastic linear system with quadratic cost (LQG), one can formulate the so-called Riccati LMI to find the optimal value function of  the LQG problem \citep{boyd1994linear,balakrishnan2003semidefinite}.  Similarly, the well known LMI formulations have not attempted to show how these results can be derived from a more general approach to stochastic optimal control, namely the inf-LP approach. 

In this work, we establish the connection between the inf-LP approach and the well-known  Riccati LMIs for LQG problems. This inf-LP in general, includes infinitely many constraints on the occupation measure.  The relaxation of these constraints to moments up to order two of the occupation measure  and taking the dual of this problem results in the well-known Riccati LMI solution approaches. Since the variables in the relaxation of primal inf-LP are discounted moments of the state and input, moment constraints or certain class of chance constraints can be naturally encoded in the inf-LP formulation. 

Our paper is organized as follows. In Section \ref{sec:problem} we review the inf-LP approach to discrete-time infinite horizon discounted stochastic control. In Section \ref{sec:approach}  we apply the approach to LQG problems. In Section \ref{sec:numerical} we provide numerical case studies.  In Section \ref{sec:conclusion} we summarize the results. 

\section{Inf-LP approach to Stochastic control}\label{sec:problem}

Consider the discrete-time  stochastic system 
\begin{align}
\label{eq:abstract_dynamics}
x_{t+1} \sim \tau(B_x|x_t,u_t),
\end{align}
where $x_t \in X$, $u_t \in U$, and $\tau(.|x,u)$ is a stochastic kernel. It assigns a probability distribution to $B_x \in \m{B}(X)$ given $x$ and $u$, where $\m{B}(X)$ is the set of Borel subsets of $X$. The  stochastic control problem is defined by
\begin{align}
\label{opt:lp_abstract}
\min \limits_{\pi \in \Pi} \; &  \mathbb{E}_{\nu_0}^\pi  \sum_{t=0}^\infty \alpha^t c_0(x_t,u_t).
\end{align}
Above, $c_0: X \times U \rightarrow \R_+$ is the running cost and $\alpha \in (0,1)$ is a discount factor, $\nu_0$ is an initial state distribution. We consider randomized policies $\pi \in  \Pi$, where $\Pi$ is the set of probability measures on $U$ given $X$. That is, for each $x \in X$, $\pi(x)$ gives a probability distribution on the input space $U$. The expectation $\mathbb{E}$ is with respect to the probability measure induced by $\nu_0$,  $\pi$ and $\tau$. 

The solution to the stochastic control problem above can be characterized as the solution of   an infinite dimensional linear program (inf-LP). To  present this inf-LP, we first define the infinite dimensional optimization spaces for the primal and dual LPs. Define the weight functions
\begin{align}
\label{eq:weights}
w(x,u) = \eps +c_0(x,u) , \; \tilde{w}(x) =\min_{u \in U} w(x,u), 
\end{align} 
where $\eps > 0$ so that the weights are bounded away from zero.  Let $\m{F}(X\times U), \m{F}(X)$ denote the space of real valued measurable functions with bounded $w, \tilde{w}$-norms, respectively. That is, for $f \in \m{F}(X\times U), \tilde{f} \in \m{F}(X)$: 
\begin{align*}
\sup_{(x,u)}\frac{|f(x,u)|}{w(x,u)} < \infty, \quad \sup_{x}\frac{|\tilde{f}(x)|}{\tilde{w}(x)} < \infty.
\end{align*}
Let $\m{M}(X\times U), \m{M}(X)$ denote the space of measures with finite $w, \tilde{w}$-variations, respectively. That is, for $\mu \in \m{M}(X\times U), \tilde{\mu} \in \m{M}(X)$:
\begin{align}
\label{eq:measure}
\int_{X\times U} w d\mu < \infty, \quad \int_{X} \tilde{w} d\tilde{\mu} < \infty.
\end{align}

Define the linear map $T: \m{M}({X\times U}) \rightarrow \m{M}(X)$ as:
\begin{align}
\label{eq:T_map}
[T\mu ](B) = \tilde{\mu}(B) - \alpha \int_{X \times U} \tau(B|x,u)\mu(dx,du), 
\end{align} 
where $\tilde{\mu}(B) := \mu(B, U)$ and $B \in \m{B}(X)$. 
Analogously, define the linear map $T^*: \m{F}(X) \rightarrow \m{F}(X \times U)$ as:
\begin{align*}
[T^*v] \;(x,u) = v(x)- \alpha \int_X \tau(dy|x,u)v(y).
\end{align*} 
Note that the second term above $\int_X \tau(dy|x,u)v(y)$, is the expectation of the function $v$ under the stochastic kernel $\tau$. One can verify that $T$ and $T^*$ are adjoint operators: 
\begin{align*}
<T^*v, \mu>_{X\times U}\, =\, <v, T \mu>_X,
\end{align*}
where the bilinear maps are given by:
\begin{align*}
<c, \mu>_{X\times U} &=  \int_{X\times U} c(x,u)\mu(dx,du), \\
<v, \nu>_X  &= \int_X v(x) \nu(dx). 
\end{align*}

In the remainder, for simplicity, we drop the subscript of $<.\;,\;.>$ since the space is  clear from the context. To formulate the inf-LP corresponding to stochastic control, we need the following standard assumptions \citep{hernandez1996discrete}. 
\begin{assumption}\
\label{asm:lp_assm}
\begin{enumerate}[(a)]
\item \label{asm:c0_compact} The cost  $c_0$ is lower semi-continuous and inf-compact, that is, for every $x \in X$, $r \in \R$, the set $\{ u \in U \; | \; c_0(x,u) \leq r\}$ is non-empty and compact. 
\item \label{asm:tau_cont} The stochastic kernel $\tau$ is weakly continuous. 
\item \label{asm:q_norm} $\sup_{X\times U} \int_X \tilde{w}(y) \tau(dy | x, u)/w(x,u) < \infty$.
\item \label{asm:nu0} $\nu_0 \in \m{M}_+(X)$. 
\end{enumerate}
\end{assumption}
Let $\m{M}_+(X \times U) \subset \m{M}(X \times U)$ denote the cone of non-negative measures. For $\nu_0 \in \m{M}_+(X)$,  the constraint on $\mu \in \m{M}(X \times U)$, denoted by $\nu_0 - T \mu = 0$  refers to 
\begin{align}
\label{const:measure}
& \nu_0(B_x) - [T\mu](B_x) =0, \quad \forall B_x \in  \m{B}(X).
\end{align}

\begin{thm} The stochastic control problem \eqref{eq:abstract_dynamics}, \eqref{opt:lp_abstract} can be equivalently formulated as the following inf-LP:
\begin{align}
\label{opt:lp_primal}\tag{P-SC}
\min_{\mu \in \m{M}_+(K)} \quad &<c_0, \mu> \\ 
\label{con:measure_inf}
\mbox{s.t. } \quad &\nu_0 - T\mu =0.
\end{align}
\end{thm}
We summarize the idea of the proof and refer the readers to \citep{hernandez1996discrete} for details.
Given a  policy $\pi \in \Pi$, one can define $\mu \in \m{M}_+(X\times U)$ as 
\begin{align}
\label{eq:discounted_ocupancy}
\mu(B_x , B_u) = \sum_{t = 0}^\infty \alpha^t \mathbb{P}^\pi_{\nu_0}\{ (x_t,u_t) \in (B_x ,B_u)\}, 
\end{align}
where $B_x \in \m{B}(X), B_u \in \m{B}(U)$. This measure corresponds to discounted probability of $(x_t,u_t)$ being in any Borel  subset of $X \times U$ and is referred to as the occupation measure.  It can be verified that the occupation measure satisfies  $\nu_0 - T\mu =0$. Furthermore, given any $\mu \in \m{M}_+(X\times U)$, there exists a policy $\varphi \in \Pi$, satisfying
\begin{align} 
\label{eq:conditional}
\mu(B_x, B_u) = \int_{B_x} \varphi(B_u|x)\tilde{\mu}(dx),
\end{align}
for all $B_x \in \m{B}(X), B_u \in \m{B}(U)$  [Proposition D.8(a) in \citep{hernandez1996discrete}]. It can be shown that the cost \eqref{opt:lp_abstract} corresponding to the policy $\varphi$ is
\begin{align}
\label{eq:cost_policy}
\mathbb{E}_{\nu_0}^\varphi  \sum_{t=0}^\infty \alpha^t c_0(x_t,u_t) = <c_0,\mu>. 
\end{align}
Putting the above results together, the problem of finding the optimal policy for \eqref{opt:lp_abstract} can be equivalently formulated as finding a measure minimizing \eqref{eq:cost_policy} subject to \eqref{con:measure_inf}.

Whereas the inf-LP above provides the optimal occupation measure and the optimal policy for the stochastic control problem, the dual of this inf-LP can be used to find the optimal value function. Furthermore, the duality gap  is zero \citep{hernandez1996discrete}. 

To define this dual inf-LP, let the constraint on $v \in \m{F}(X)$, denoted by $c_0-T^*v \geq 0$ refer to
\begin{align}
\label{const:valuefunction}
& c_0(x,u)- [T^*v](x,u) \geq 0, \; \forall (x,u) \in X\times U. 
\end{align}
The dual inf-LP is given by:
\begin{align}
\label{opt:lp_dual}\tag{D-SC}
\max\limits_{v \in \m{F}(X)} \quad &  < v, \nu_0>  \\
\label{con:constraint_u_inf}
\mbox{s.t. } \quad &  c_0 -T^*v\geq 0.
\end{align}
\textbf{Remark.}  Constraint \eqref{con:constraint_u_inf} is   the Bellman inequality. In particular, based on the Bellman principle of optimality, a function $v^*$ is the optimal value function of the stochastic control if and only if $ c_0 -T^*v= 0$. Thus, the optimizer of the above inf-LP satisfies the Bellman equality.

\section{Inf-LP Approach to LQG Problems}\label{sec:approach}
Consider the linear system as a specialization of \eqref{eq:abstract_dynamics}:
\begin{align}
\label{eq:linear_dynamics}
x_{t+1} = Ax_t + B u_t + \omega_t, 
\end{align}
where $x \in \R^n$, $u \in \R^m$, $\omega \in \R^n$,  $\omega_t$, are independent identically distributed Gaussian random variables and for all $t$, $E\{\omega_t\} = 0$ and $E\{\omega_t \omega_t^T\} = W$.  The initial state is independent of the stochastic noise and has a distribution $\nu_0$, with mean $E\{x_0\} = m_0$ and covariance $E\{x_0 x_0^T\} = \Sigma_0$. The discounted linear quadratic Gaussian (LQG) problem is formulated as:
\begin{align}
\label{eq:lqg_cost}
\min \limits_{\pi \in \Pi} \;  &\mathbb{E}_{\nu_0}^\pi  \sum_{t=0}^\infty \alpha^t (x^T_t Q_0 x_t  + u^T_t  R_0 u_t ).
\end{align}
We assume the pair $(A,B)$ is controllable and the pair $(A,C)$ is observable, where $Q_0= C^TC$. Denote by $\m{S}^n$, $\m{S}^n_+$ and $\m{S}^n_{++}$  the set of $n \times n$ symmetric, symmetric positive semidefinte and symmetric positive definite matrices, respectively.  We assume $Q_0 \in \m{S}^n_+$ and $R_0 \in \m{S}^m_{++}$.

To apply the inf-LP approach to this problem, first we verify Assumption \eqref{asm:lp_assm} as follows. The weight functions \eqref{eq:weights} in the LQG problem are  $w(x,u) = \eps + x^T Q_0  x + u^T  R_0  u$, $\tilde{w}(x) = \eps + x^T Q_0x$. Thus, $\m{F}(X\times U)$ and $\m{F}(X)$ are spaces of functions over $X\times U = \R^{n\times m}, X = \R^n$ respectively, that do not grow faster than quadratic functions. Furthermore, by definition \eqref{eq:measure}, $\m{M}(\R^{n\times m}), \m{M}(\R^n)$ are sets of measures  that have bounded variance. From $Q_0 \in \m{S}_+, R_0 \in \m{S}_{++}$, part (a) of Assumptions \eqref{asm:lp_assm} holds.  The stochastic kernel $\tau$ is Gaussian, which is continuous and has finite variance, satisfying part (b). Part (c) holds since 
\begin{align*}
& \int_X w_0(y) \tau(dy | x, u) = \epsilon + x^T A^T Q_0Ax \\
&+ 2x^TA^TQ_0Bu + u^TB^TQBu+ \text{Tr}(Q_0W) \in \m{F}(\R^{n\times m}).
\end{align*}
Finally, part (d) holds due to finite variance of the initial  state distribution $\nu_0$. 

\subsection{Primal and Dual SDPs}
We consider a relaxation of the inf-LP \eqref{opt:lp_primal}, resulting in an equivalent restriction of \eqref{opt:lp_dual}, to obtain a tractable formulation of these inf-LPs. These formulations are then connected to the well-known Riccati LMIs to solve the LQG \citep{balakrishnan2003semidefinite}. 

First, consider the constraint $\nu_0-T\mu = 0$ in \eqref{opt:lp_primal}. It can be verified that this is equivalent to $ <v,\nu_0-T\mu> = 0, \; \forall v \in {\m{F}}(X)$. We relax this constraint by restricting $\m{F}$ to a subset $\hat{\m{F}}(X)$. In particular, define
\begin{align}
\label{eq:quadratic_functions}
\hat{\m{F}}(X)= \{ v \in \m{F}(X) \; | \; v(x) = x^T P x + q^Tx + r \},
\end{align}
$P \in \m{S}^n, q \in \R^n, r \in \R$. Since $v \in \m{F}(X)$ are quadratic, the infinitely many constraint \eqref{con:measure_inf} on the measure $\mu$ are relaxed to a set of finite constraints on its moments of order up to two. These constraints will be derived as follows. 

Introduce moments of measure $\mu$ as:  
\begin{align*}
m & = \int_{X \times U} \mu(dx, du) = \mu(X\times U) \in \R_+, \\
m_x & = \int_X x \mu(dx,U) \in \R^n, \\
Z_{xx} & = \int_X x x^T \mu(dx,U) \in \m{S}^n_+. 
\end{align*}
Similarly, $m_u, Z_{xu}, Z_{uu}$ are defined. For any $S \in \R^{n\times n}$, $\text{Tr}(S)$ denotes trace of the matrix $S$.

\begin{prop}[Primal SDP for LQG]
\label{lem:primal_sdp}
By constraining $\m{F}(X)$ to $\m{\hat{F}}(X)$, we obtain a relaxation of \eqref{opt:lp_primal} as
\begin{align}
\label{sdp:primal}\tag{P-LQG}
& \min \quad \text{Tr}(Q_0 Z_{xx}) +\text{Tr}( R_0 Z_{uu})  \\ 
\nonumber
\mbox{s.t. } \quad  & \Sigma_0 + m_0 m_0^T - Z_{xx} + \alpha A Z_{xx}A^T 
 + A Z_{xu}B^T \\
\label{con:quad_measure}\tag{C2}
&+ B Z_{xu}^T A^T + B Z_{uu} B^T +m W ) = 0_{n\times n}\\
\label{con:mean}\tag{C1}
\quad  & m_0 - m_x + \alpha(A m_x + Bm_u) = 0_{n\times 1},\\
\label{con:total}\tag{C0}
\quad  & \alpha m -m + 1 =0,\\
\label{con:quad_measure2}\tag{Cp}
\quad  &Z:= \left[  \begin{array}{ccc} m & m_x^T & m_u^T \\ m_x & Z_{xx} & Z_{xu} \\ m_u & Z_{xu}^T & Z_{uu}\end{array} \right] \succeq 0,
\end{align}
over the variables $(m,m_x, m_u, {Z_{xx}, Z_{xu}, Z_{uu}})$. 
\end{prop}
\begin{pf}
Based on the definition of the moments, the cost functions in Problem \eqref{eq:lqg_cost} can be expressed as:
\begin{align*}
< c_0, \mu> = \text{Tr}(Q_0 Z_{xx}) + \text{Tr}(R_0 Z_{uu}).
\end{align*}
Next, expanding $ <v,\nu_0-T\mu> = 0$ we obtain the first term as
\begin{align*}
 <v,\nu_0 >&= <x^TPx + q^T x + r, \nu_0> \\
&=  \text{Tr}(P \Sigma_0) + m_0^TPm_0 + q^T m_0 + r.
\end{align*}
Using definition of $T$ in \eqref{eq:T_map}, we expand $<v, -T\mu>$. The first term is:
\begin{align*}
&<v, \tilde{\mu}>\,= \,<x^TPx + q^T x + r, -\tilde{\mu}> = \\
&-  \int_X v(x) \mu(dx,U) = - \text{Tr}(PZ_{xx}) - q^T m_x - m r.
\end{align*}
The second term is obtained as: 
\begin{align*}
&\alpha \times \big(\text{Tr}(P (A Z_{xx}A^T + AZ_{xu}B^T + B Z_{xu}^T A^T  + B Z_{uu}B^T) \big)  \\
&+    q^T(Am_x + B m_u) + m \big(\text{Tr}(PW) + r\big).
\end{align*}
In the above, we used the fact that the measure $\tau(dx|y,u)$ has mean $Ay + Bu$ and covariance $W$. 
Each of the terms in the constraint expansion above are linear in the variables $P, q, r$. Since $<v, \nu_0-T\mu>=0$ must hold for all $v \in \hat{\m{F}}(X)$, that is, for all $P \in S^{n}, q \in \R^n, r \in \R$, the corresponding coefficients of these variables need to equal zero. From this, we obtain the set of affine constraints,  \eqref{con:total}, \eqref{con:mean}, \eqref{con:quad_measure}. Constraint \eqref{con:quad_measure2} holds since $Z$ is moment of a positive measure $\mu$ \citep{lasserre2009moments}. \qed
\end{pf} 

Similarly, we can obtain the dual SDP as follows. 
\begin{prop}[Dual SDP for LQG]
\label{lem:dual_sdp}
By constraining $\m{\hat{F}}(X)$ to $\m{F}(X)$, we obtain a restriction of \eqref{opt:lp_dual}  as follows:  
\begin{align}
\label{sdp:dual}\tag{D-LQG}
& \min \; \text{Tr}(P\Sigma_0) + \text{Tr}(Pm_0 m_0^T) + q^Tm_0 +  r  \\ 
\nonumber
& \mbox{s.t. }  \quad \left[ \begin{array}{ccc} s_0 & s_1^T & s_2^T  \\ s_1 & S_{11} & S_{12} \\ s_2 & S_{12}^T &  S_{22} \end{array} \right] \succeq 0, 
\end{align}
with optimization variables, $P, q, r$ and 
\begin{align*}
& s_0= r(\alpha-1) + \alpha \text{Tr}(PW), \\
& s_1 =  \frac{1}{2} (-I + \alpha A^T)q, \; s_2 = \frac{\alpha}{2}  B^Tq,  \\
&\small{\left[ \begin{array}{cc}S_{11} & S_{12} \\   S_{12}^T & S_{22} \end{array} \right] = \left[ \begin{array}{cc} \alpha A^T P A-P + Q_0  & \alpha A^T P B \\   \alpha B^T P A & R_0 + \alpha B^T P B\end{array} \right]. }
\end{align*}
Furthermore, this SDP is the dual of \eqref{sdp:primal}.
\end{prop}
\begin{pf}
The term $<v, \nu_0>$ in the cost function was discussed in Proof of Proposition \eqref{lem:primal_sdp}. For $v \in \hat{\m{F}}(X)$,  Constraint \eqref{const:valuefunction} becomes
\begin{align*}
& x^T (-P + \alpha A^T P A + Q_0)x+ 2 \alpha  x^T A^T P B u  \\
+& u^T (\alpha  B^T P B + R_0)u + q^T(-I + \alpha A)x \\
+ & \alpha q^T B u + r(\alpha-1) + \alpha \text{Tr}(PW) \geq 0, \; \forall (x,u) \in X \times U. 
\end{align*}
An equivalent way of writing the above constraint is:
\begin{align*}
 & \left[ \begin{array}{c}  1 \\ x \\ a \end{array} \right]^T 
\left[ \begin{array}{ccc} s_0 & s_1 & s_2  \\ s_1^T & S_{11} & S_{12} \\ s_2^T & S_{12}^T &  S_{22} \end{array} \right] 
\left[ \begin{array}{c}  1 \\ x \\ a \end{array} \right]  \succeq 0,
\end{align*}
which leads to the constraint in \eqref{sdp:dual}.  Using SDP duality \citep{vandenberghe1996semidefinite}, it can also be verified that \eqref{sdp:dual} is dual of \eqref{sdp:primal}.\qed
\end{pf}
\textbf{Remark.} The primal and dual SDPs above are a generalization of the existing results in literature  due to the  additional terms arising from zero and first order moments of $\mu$. If $m_0 = 0$  from controllability of the pair $(A,B)$ we have $m_x = 0$, $m_u = 0$. Thus, removing the corresponding dual variable $q \in \R^n$, we can obtain that $r  = \frac{\alpha}{1-\alpha}\text{Tr}(PW)$. This leads to the standard results in \citep{willems1971least,boyd1994linear,balakrishnan2003semidefinite}:
\begin{align}
\label{sdp:dual2}\tag{D-LQG0}
& \min \; \text{Tr}(P\Sigma_0) + \frac{\alpha}{1-\alpha}\text{Tr}(PW)   \\ 
\nonumber
& \mbox{s.t. }  \quad S:= \left[ \begin{array}{cc}  A^T P A-P + Q_0& \alpha A^T P B \\  \alpha B^T P  & R_0 + \alpha B^T P B\end{array} \right] \succeq 0. 
\end{align}
In the rest of the paper, we consider $m_0=0$ and thus, we work with \eqref{sdp:dual2} and its dual. 

If the SDP \eqref{sdp:dual2} and its  dual have non-empty optimal sets, the complementary slackness condition holds \citep{vandenberghe1996semidefinite}:
\begin{align*}
\nonumber
&Z^* S^* = 0 \iff \left[ \begin{array}{cc} Z_{xx} & Z_{xu} \\   Z_{xu}^T & Z_{uu} \end{array} \right] \times \\
 & \left[ \begin{array}{cc} -P  + Q_0 + \alpha A^T P A & \alpha A^T P B \\   \alpha B^T P A & R_0 + \alpha B^T P B\end{array} \right] = 0, 
\end{align*}
where we dropped $^*$ from individual terms above. Expanding above equality, we obtain that
\begin{align}
\nonumber
0 = \;& -P  + Q + \alpha A^T P A +\\
\label{eq:p_riccati}
&  \alpha A^T P B (R + \alpha B^T P B)^{-1}\alpha B^T P A, \\
\label{eq:lqg_gain}
 Z_{xu}^T Z_{xx}^{-1}  = \;& (R + \alpha B^T P B)^{-1}\alpha B^T P A, \\
\label{eq:covariance_0}
Z_{uu} = \;& Z_{xu}^T Z_{xx}^{-1}Z_{xu}.
\end{align}
%
Equation \eqref{eq:p_riccati} is the algebraic Riccati equation of the infinite horizon discounted LQG problem. Equation \eqref{eq:lqg_gain} provides the optimal controller gain, $K =  Z_{xu}^T Z_{xx}^{-1}$. 

\textbf{Remark.} An alternative derivation of the optimal policy is provided by considering the occupation measure $\mu$ in the inf-LP. Since $ \frac{1}{m}\mu({X \times U}) = (1-\alpha)  \mu({X \times U}) $ is a Gaussian  measure (alternatively, by considering only the knowledge of the first and second order moments of this measure),  the conditional measure $\varphi(u|x)$ \eqref{eq:conditional} can be obtained as 
\begin{align*}
\varphi(u|x) \sim \m{N}(m_{u|x},Z_{u|x}), 
\end{align*}
with the conditional mean $m_{u|x} = m_u +  Z^T_{xu} Z_{xx}^{-1} x$ and covariance $Z_{u|x} =  \frac{1}{m} (Z_{uu} - Z^T_{xu}Z_{xx}^{-1} Z_{xu})$.  By complementary slackness of \eqref{eq:covariance_0}, the covariance of this measure is zero and thus, the optimal policy predicted by inf-LP \eqref{opt:lp_primal} is deterministic and is equal to $\varphi(x) = Z^T_{xu} Z_{xx}^{-1} x$.

\subsection{Constrained LQG}\label{sec:clqg}
One of the advantages of the inf-LP  \eqref{opt:lp_primal} is that  constraints of the form $\mathbb{E}_{\nu_0}^\pi  \sum_{t=0}^\infty \alpha^t c_i(x_t,u_t)  \leq \beta_i$, $i = 1, 2, \dots, N$, can readily be incorporated. In particular, from the definition of the occupation measure  \eqref{eq:discounted_ocupancy}: $$\mathbb{E}_{\nu_0}^\pi  \sum_{t=0}^\infty \alpha^t c_i(x_t,u_t)  \leq \beta_i \iff <c_i, \mu> \; \leq \; \beta_i.$$ Such constraints  correspond to multi-objective stochastic control, where $c_i$, for $i=1, \dots, N$,  denotes a set of  additional objectives and $\beta_i \in \R_+$ are desired bounds.  

For the LQG problem, let $c_1(x_t,u_t) = x_t^TQ_1x_t + u_t^T R_1 u_t$. Then, constraint $<c_1, \mu> \; \leq\; \beta_1$  is equivalent to $\text{Tr}(Q_1Z_{xx}) + \text{Tr}(R_1Z_{uu}) \leq \beta_1$ in the primal SDP \eqref{sdp:primal}. From our derivation of \eqref{sdp:dual2}, it can be verified that the dual SDP with the additional constraint is
\begin{align}
\label{sdp:cdual}\tag{C-LQG}
& \min \; \text{Tr}(P\Sigma_0)   + \frac{\alpha}{1-\alpha}\text{Tr}(PW) -  \gam \beta \\
\nonumber
& \mbox{s.t. }  \quad \left[ \begin{array}{cc} -P  + Q + \alpha A^T P A & \alpha A^T P B \\   \alpha B^T P A & R + \alpha B^T P B\end{array} \right]\succeq 0, 
\end{align}
where $Q = Q_0 + \beta Q_1$ and $R = R_0 + \beta R_1$ and the optimization variables are $P$ and the dual multiplier of the constraint $\gam > 0$. This is consistent with alternative derivations in multi-criterion LQG in \citep{boyd1994linear}.  

Due to $Z_{xx}$ and $Z_{uu}$ corresponding to the second order discounted moments of the occupation  measure,   constraints $\text{Tr}(Q_1Z_{xx}) \leq \beta_1$ can also be used to pose chance constraints on the state of the form
\begin{align}
\label{eq:chanc_constraint}
\sum_{t=0}^\infty \alpha^t\mathbb{P}_{\nu_0}^\pi  ( |g^T x_t |< h) \geq 1-\eps.
\end{align} 
\begin{prop}
Given a policy $\pi$, let $\mu$ be the resulting occupation measure. Let $Q_1 = gg^T$, and $\beta_1 = \eps h^2$
\begin{align*}
\text{Tr}(Q_1Z_{xx}) \leq \beta_1 \Rightarrow \sum_{t=0}^\infty \alpha^t P^\pi_{\nu_0} ( |g^T x_t |\leq h) \geq 1-\eps. 
\end{align*}
Furthermore, for the case in which $\pi(x) = Kx$ and in the LQG setting,  letting $\beta_1 =  \frac{h^2}{ 2 (\text{erf}^{-1}(1-\eps/m)^2) }$,where \text{erf} denotes the error function, the constraint above is also sufficient for the chance constraint. \end{prop}
\begin{pf}
Define $X_h = \{ x \; : \; |g^T x |\geq h \}$, so that $\mu(X_h, U) = \sum_{t=0}^\infty \alpha^t \mathbb{P}_{\nu_0}^\pi \big( x_t \in X_h \big)$. Then, the chance constraint \eqref{eq:chanc_constraint} can be written as $\mu(X_h,U) \leq \epsilon$. 
Now
\begin{align*}
0 \leq h^2\mb{1}_{X_h} \leq (g^T x)^2 \mb{1}_{X_h}
\end{align*}
Taking integral with respect to $\mu$ we obtain that
\begin{align*}
h^2\mu({X_h}, U) &\leq \int_X (g^T x)^2 \mb{1}_{X_h}\mu(dx,U) \leq g^TZ_{xx}g.
\end{align*}
It follows that $\frac{g^TZ_{xx}g}{h^2} < \eps \Rightarrow \mu(X_h,U) < \eps$.

If $\pi$ is linear, the resulting occupation measure $\mu$ is a scaled (by factor $m = \frac{1}{1-\alpha})$ Gaussian measure.  Since the cumulative distribution function of the Gaussian distribution is invertible and is given through the \text{erf} function, we have
\begin{align*}
\mu(X_h, U) < \eps & \iff 1 + \text{erf}(\frac{-h}{\sigma \sqrt{2}}) \leq \frac{\epsilon}{m} \\
 &  \iff \frac{h^2}{ 2 (\text{erf}^{-1}(1-\eps/m)^2) } \geq g^TZ_{xx}g. \qed
\end{align*}
\end{pf}
\textbf{Remark.} The above chance constraints were also considered in infinite horizon average cost LQG problems \citep{schildbach2015linear}.  The authors derived  analogous SDPs for the average cost criterion. Their approach was through considering the steady-state second order moments of the closed loop linear system, corresponding to a linear policy. As shown here,  this approach is equivalent to  the relaxation of the primal infinite dimensional LP from which the occupation measure and the corresponding second order moments were derived. 


\section{Numerical Case Studies} \label{sec:numerical}
Our goal is to use the constrained LQG formulation of the previous section to study effects of multi-objective and chance constraints on the LQG problem. To this end, we consider two linear systems, each with a nominal objective, and closed-loop infinite horizon second order moment constraints on states or inputs. In both cases, we solved SDP \eqref{sdp:primal} using the parser CVX \citep{grant2008cvx} with the solver SDPT3 \citep{toh1999sdpt3}.

Our first example is  a second order system to study effects of chance constraints  \eqref{eq:chanc_constraint}. In our second example, we consider a model for a miniature coaxial helicopter linearized around a hover maneuver, with a nominal objective of minimizing deviations from hover. Our secondary objective is to minimize control energy. 

\subsection{Two-state system with state constraints}
The system dynamics parameters are 
\begin{align*}
A = \left[ \begin{array}{cc}
   1 \;& 0.1 \\
   0 \;& 1 \\
  \end{array} \right], \; 
B = \left[ \begin{array}{c}
     0 \\
     1 \\
    \end{array} \right],  \; 
W = \left[ \begin{array}{cc}
         0.1 & 0 \\
          0 & 0.1\\
        \end{array} \right].
\end{align*}
The primary and secondary objective parameters are
\begin{align*}
Q_0 = \left[ \begin{array}{cc}
   1 & 0 \\
   0 & 1 \\
  \end{array} \right], \; 
Q_1 = \left[ \begin{array}{cc}
   0 & 0 \\
  0 & 1 \\
  \end{array} \right], \; 
 R_0 = 1, R_1 = 0.
\end{align*}
The discount factor is $\alpha = 0.99 $ and $x_0 \sim \m{N}(m_0, \Sigma_0)$, with mean $m_0 = [-0.46, 0.58]^T$ and covariance $\Sigma_0 = I$. The second objective, $\text{Tr}(Q_1Z_{xx})$ is constrained to be less than a parameter $\beta$. As such, we require $\sum_{t=0}^\infty \alpha^t P^\pi_{\nu_0} ( [0,\,1] x_t \leq h) \geq 1-\eps$. This has the interpretation of a soft constraint for the second state to remain close to zero. 

\begin{figure} \label{fig:tradeoff1}
\begin{center} 
\resizebox{1.09\linewidth}{!}{\includegraphics{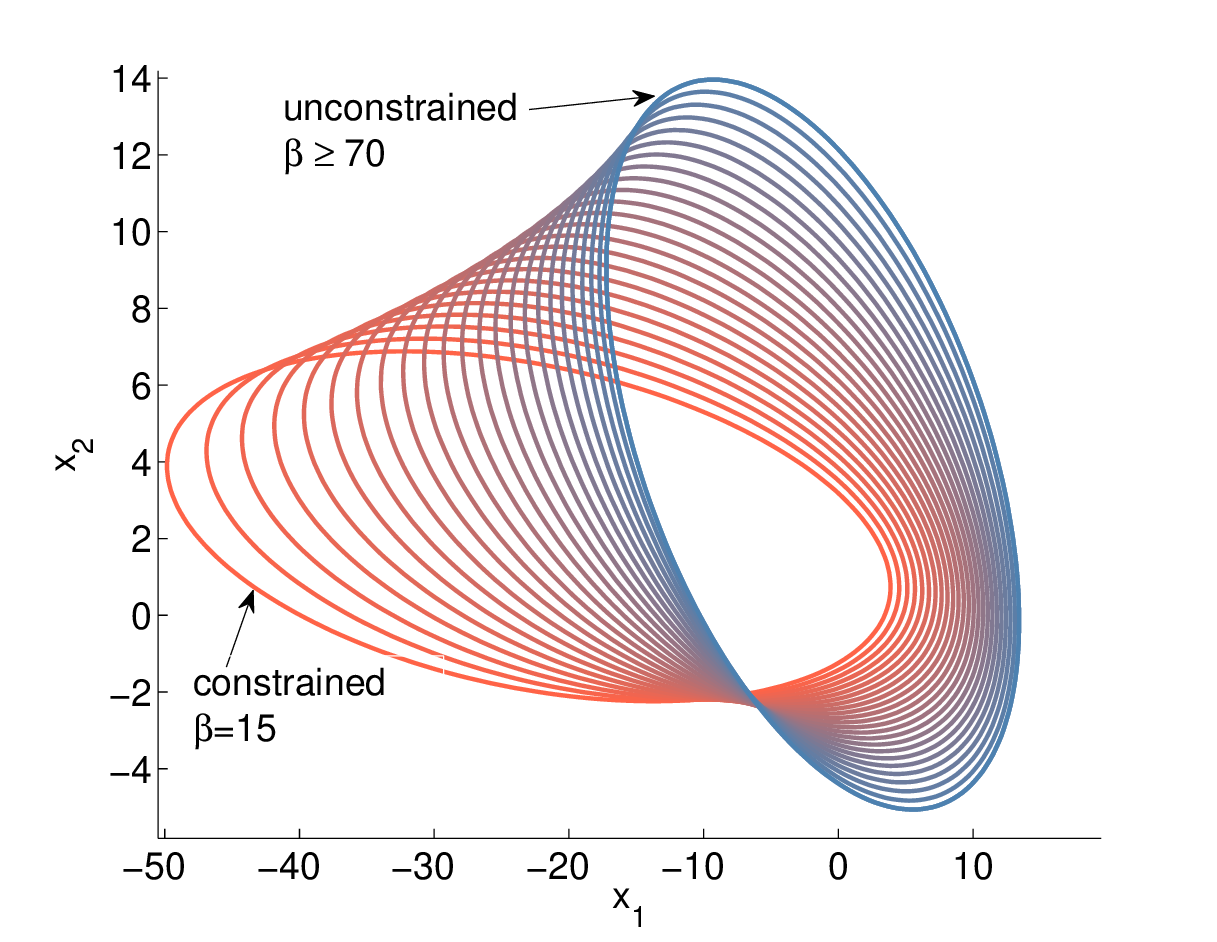}}
\caption{Discounted state occupancy ellipses $\{x \in \mathbf{R}^2 \mid (x - m_x^*)^T Z_{xx}^{*-1} (x - m_x^*) = 1 \}$ as the constraint associated with the secondary cost is tightened from unconstrained (blue) to $\beta = 15$ (red).}
\end{center}
\end{figure}

We vary $\beta$ between the value of the second objective achieved using the optimal policy for first objective without the constraint and a lower tighter value that forces the constraint to be active. Figure \eqref{fig:tradeoff1} shows the optimal state covariance $Z^*_{xx}$. It is seen that  as the constraint is tightened, the discounted occupancy ellipse changes to one with less variation in the second state, due to definition of $Q_1$, but more variation in the first state. The cost is 378 without the constraint and 981 with $\beta = 15$.  

\subsection{Miniature coaxial helicopter}
We now consider a simplified eight-state model of a miniature two-rotor coaxial helicopter linearized around a hover maneuver based on \citep{kunz2013fast,summers2013approximate}. The states of the system are the three-dimensional position and heading deviations from a desired hover pose in an inertial reference frame and the associated velocities in a body reference frame. There are four inputs: pitch, roll, thrust, and yaw, used for forward flight, sideways flight, vertical flight, and heading change, respectively. Pitch and roll are actuated with a swashplate mechanism connected to two servos. Thrust is actuated by the rotational speed of the rotor motors, and yaw is actuated by a rotational speed difference of the rotor motors. The pitch and roll angles and velocities are neglected in the model, and the pitch and roll inputs are assumed to act directly on the lateral position states. 

The dynamics are discretized in time by Euler integration with sampling time $t_s$. The parameters of the the discrete-time system dynamics are 
\begin{align*}
&A = \left[ \begin{array}{cc}
   I_4 & t_s I_4   \\
   0 & I_4   + t_s \text{diag} ([k_x, k_y, k_z, k_\psi])  \\
  \end{array} \right], \;  \\
&B = \left[ \begin{array}{c}
     0_4 \\
     t_s \text{diag} ([b_x, b_y, b_z, b_\psi])  \\
    \end{array} \right],  \; 
W = \left[ \begin{array}{cc}
         0_4  & 0 \\
          0 & 0.1I_4  \\
        \end{array} \right],
\end{align*}
where $[k_x,k_y,k_z,k_\psi ] = [-0.5,\, -0.5,\, 0,\, -5]$ represent fuselage drag parameters and $[b_x,b_y,b_z,b_\psi] = [2.0,2.1,11,18]$ represent inertial parameters mapping actuator influence to state derivatives. The values of the parameters are taken from \citep{kunz2013fast} and are based on a grey-box system identification with experimental data.  

We consider a scenario in which the deviations from the desired hover pose are to be minimized, subject to an infinite-horizon closed-loop constraint on the expected discounted control energy. The trade-off between trajectory optimization and energy cost minimization is a classical control tradeoff. It is encoded in our framework with the primary and secondary cost parameters
\begin{align*}
Q_0 = I_8 , \; 
Q_1 = 0_8 , \; 
R_0 = 0_4 , 
R_1 = I_4 .
\end{align*}
The discount factor is $\alpha = 0.99$ and the threshold is $\beta = 50$. 
The constraint on control energy is achieved with the constrained LQG formulation by solving \eqref{sdp:primal}. This constraint is satisfied at a price of poorer regulation compared to the unconstrained case as illustrated by the closed-loop system performance in Fig. 2 .
\begin{figure} \label{fig:tradeoff2}
\begin{center} 
\resizebox{1.03\linewidth}{!}{\includegraphics{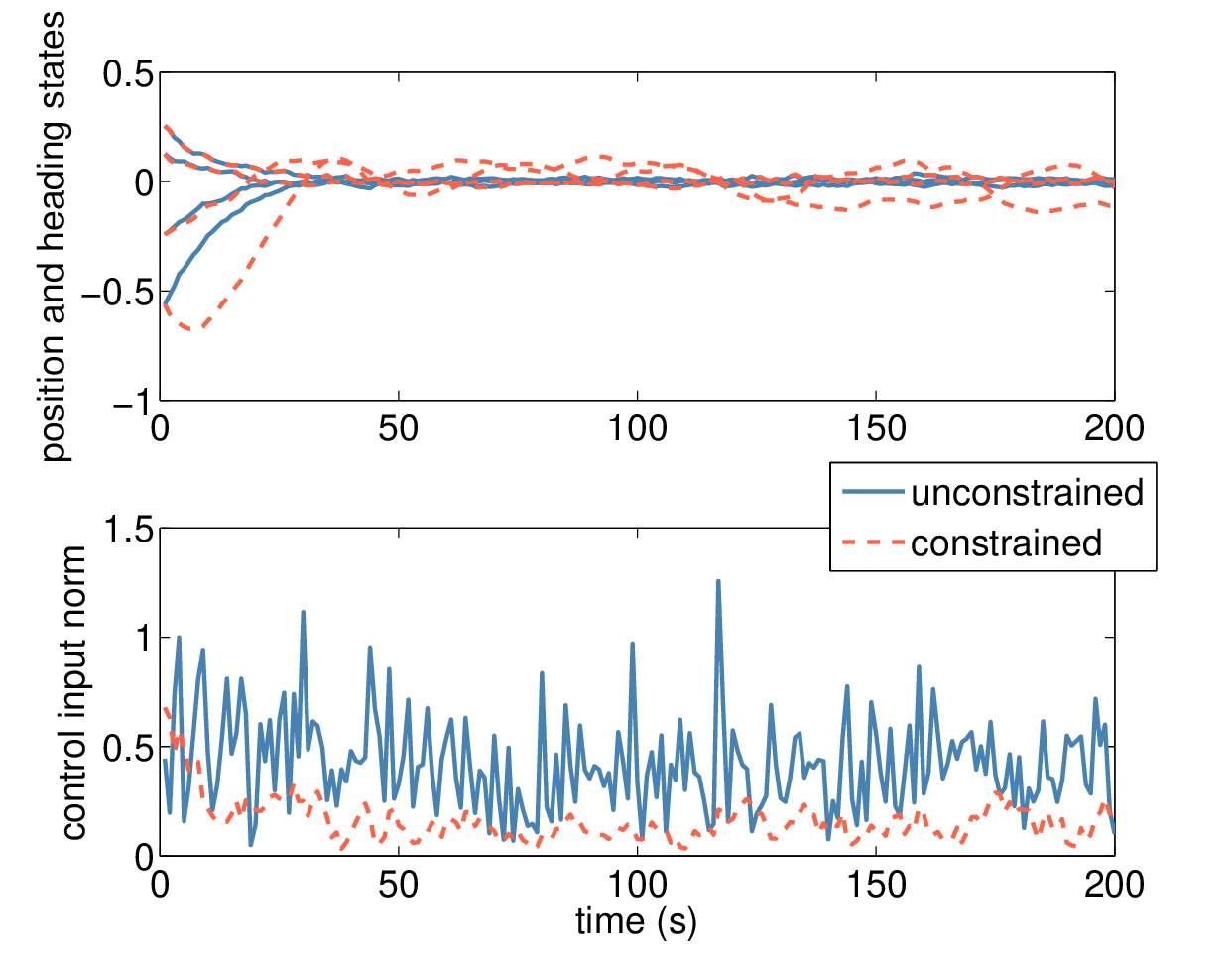}}
\caption{Classical state regulation and control energy tradeoff achieved with the constrained LQG formulation. The plots show a realization of the position and heading state and control norm evolution with and without the  constraint on the closed-loop infinite-horizon expected discounted control energy.}
\end{center}
\end{figure}

\section{Conclusions}\label{sec:conclusion}
We established a link between  the semidefinite programs (SDPs) for solving the LQG problem and the infinite dimensional linear programming (inf-LP) approach to stochastic control. The inf-LP approach is an equivalent alternative to the Dynamic Programming  principle of optimality. While the inf-LP formulation has been known since 1950's, its computational aspects and connections with existing control theoretic results  have not been fully explored. We showed that the LMI derived from the occupation measure formulation of the inf-LP corresponds to the dual of the well-known Riccati LMI. Furthermore, given second order moments of the occupation measure, we showed that multi objective and chance constraints have a natural interpretation in this framework, and these formulations coincide with alternative approaches to derive the results.  We illustrated the constrained LQG problem with two numerical case studies. 

Extensions of this work to continuous-time LQG, and average cost LQG problems are straightforward and will complete the picture. While a rich theory of approximate dynamic programming (ADP) exists, it would be interesting to enrich the approximation procedures, through advanced optimization techniques for solving the inf-LPs corresponding to stochastic control. There has been recent promising steps towards this objective \citep{sutter2014approximation,adp_reach_nikos}. It will be interesting to further apply these techniques to  large-scale and constrained stochastic control problems. 

\bibliography{references_maryamka}

\end{document}